\documentclass[12pt]{amsart}

\usepackage{amssymb,latexsym}

\numberwithin{equation}{section}

\def\cal{\mathcal}

\def\Bbb{\mathbb}

\def\C{{\Bbb C}}

\def\N{{\Bbb N}}
\def\R{{\Bbb R}}

\def\PSet{\mbox{\rm I\kern-.22em P}}

\def\al{{\alpha}}
\def\S{{\cal S}}

\def\j{\jmath}
\def\D{\mathcal{D}}
\def\la{\lambda}

\def\Im{{\rm Im}\,}
\def\Re{{\rm Re}\,}

\def\({( \hspace{-0.335em}(}
\def\){) \hspace{-0.335em})}
\def\supp{{\rm supp\,}}

\def\sinh{{\rm sh}}
\def\cosh{{\rm ch}}
\def\arccosh{{\rm arcch}}

\def\dst{\displaystyle}
\def\gam{\gamma}
\def\vphi{\varphi}
\def\ve{\varepsilon}

\def\fu2{\frac{n}{2}}

\def\l({\left(}
\def\r){\right)}
\def\be{\begin{enumerate}}
\def\ee{\end{enumerate}}

\newtheorem{theorem}{Theorem}[section]

\newtheorem{lemma}[theorem]{Lemma}
\newtheorem{corollary}[theorem]{Corollary}
\newtheorem{proposition}[theorem]{Proposition}

\newtheorem{thm}{Theorem}[section]

\newtheorem{cor}[thm]{Corollary}
\newtheorem{remark}[thm]{Remark}

\begin{document}

\title{Wave equation and multiplier estimates on\\ $ax+b$ groups}
\author[D. M\"uller]{Detlef M\"uller}
\address{Mathematisches Seminar, C.A.-Universit\"at Kiel,
Ludewig-Meyn-Strasse 4, D-24098 Kiel, Germany}
\email{{\tt mueller@math.uni-kiel.de}}
\author[C. Thiele]{Christoph Thiele}
\address{Department of Mathematics, UCLA, Los Angeles CA 90095-1555}
\email{{\tt thiele@math.ucla.edu}}
\thanks{Part of this project was done while the authors were members of the
Erwin Schr\"odinger Institute in Vienna.
 The second author was supported by NSF grants DMS 9970469 and DMS 9985572 and
by a Sloan Fellowship. }

\begin{abstract}Let $L$ be the distinguished Laplacian on certain
semidirect products  of $\R$ by $\R^n$ which are of $ax+b$-type.
We prove pointwise estimates for the convolution kernels of
spectrally localized wave operators of the form
$$e^{it\sqrt{L}} \psi(\sqrt{L}/\lambda)$$
for arbitrary time $t$ and arbitrary $\lambda>0$, where $\psi$ is a smooth
bump function supported in $[-2,2]$ if $\lambda\le 1$ and
supported in $[1,2]$ if $\lambda\ge 1$. As corollary, we reprove
a basic multiplier estimate from \cite{hebisch-steger}
for this particular class of groups, and derive
Sobolev estimates for solutions to the wave equation associated
to $L$. There appears no dispersive effect with respect to the $L^\infty$ -
norms
for large times in our estimates, so that it seems unlikely that
non-trivial Strichartz type
estimates hold .

\end{abstract}

\maketitle

\section[]{Introduction}

We denote by $G$ the semi-direct product $G=\R\ltimes\R^n,$ endowed with
the group law
\[(x,y)(x',y')=(x+x',y+ e^x y').\] This subgroup
of the  affine group of $\R^n$ is a solvable Lie group
with exponential volume growth. We call $G$ an $ax+b$ group.

A basis for the Lie algebra of left-invariant vector fields is given by
\begin{equation}\label{basis}
X= \partial_x ,\ Y_1= e^x \partial_{y_1},\ \dots, \ Y_n= e^x \partial _{y_n}.
\end{equation}
We define the distinguished left-invariant Laplacian as the second order
differential operator
\begin{equation}\label{distlap}
L=- X^2 - \sum_{j=1}^n Y_j^2.
\end{equation}
A right-invariant Haar measure on $G$ is given by
$$dg=dxdy_1\dots dy_n.$$
We will use this right-invariant measure in notions such as $L^p(G)$.
An operator on
function spaces on $G$ is given by a right convolution kernel $k$ if
\begin{equation}\label{rightconv}
Tf(g')=f\ast k(g')=\int f({g}^{-1})k(gg')\, dg, \quad \text{ for} \  f\in\D
(G).
\end{equation}

The distinguished Laplacian $L$ has a self-adjoint extension in $L^2(G)$
(\cite{nelson}), and thus we can use spectral calculus to define the operators
\begin{equation}\label{smoothed}
e^{it\sqrt{L}}m(L),
\end{equation}
where the multiplier $m$ will lie in a suitable symbol class.
The main purpose of this article is to prove Theorem \ref{ktsatz},
which states uniform (w.r. to $\lambda$
and $t$) pointwise estimates for the convolution kernels of spectrally
localized  multiplier operators of the form
\begin{equation}\label{epsi}
e^{it\sqrt{L}}\psi(\sqrt{L}/\la),
\end{equation}
for arbitrary time $t\in \R$ and arbitrary $\lambda>0,$ where $\psi$ is a
bump function
supported  in $[-2,2],$ if $\la\le 1,$ and in $[1,2],$ if $\la\ge 1$.
We can use these  estimates to give a new proof of the
basic multiplier estimate used in \cite{hebisch-steger} (see Theorem 6.1 of
\cite{hebisch-steger}),
which is
based   entirely on the wave equation, at least for the class of $ax+b$ groups
under consideration.

As a corollary of our main theorem, we shall also deduce Sobolev estimates
for solutions
to the wave equation associated to $L$.

Our estimates are mainly controlled by the left-invariant Riemannian
distance
\begin{equation}\label{defr}
R=R(x,y):=\arccosh(\cosh(x)+\frac 12 \|y\|^2 e^{-x}) ,\end{equation}
of a point $(x,y)$ to the identity element on $G$, where
$y=(y_1,\dots,y_n)$,
$\|y\|^2=\sum_{j=1}^n y_j^2$ and
$\arccosh$ is  understood to take $[1,\infty)$  to $[0,\infty).$

We remark that, for  $n=2$, our main theorem could also be deduced from
 a transfer principle of Hebisch \cite{hebisch}. In that special situation, the
group $G$ is of the form $AN$ where $KAN$ is the Iwasawa decomposition of the
complex Lie group $SL_2(\C),$ and  Hebisch  introduces a mapping from radial
functions $f$ on $\R^3$ to functions on the group $G,$ given by
$$Tf(x,y)= C e^{-x} \frac R{\sinh(R)}f(R),$$
 which preserves the $L^1$ norm and commutes with convolution and with
application of the corresponding Laplacians.
 This allows to deduce our main Theorem \ref{ktsatz} from
the analoguous theorem on $\R^3$. However, this transfer principle is somewhat
misleading, since it would suggest for higher dimensions  estimates different
from the ones which actually hold on $G.$

We also remark that our results should extend to the distinguished Laplacians
which arise from the Laplace-Beltrami operators on rank one symmetric spaces
of the non-compact type by means of conjugation with the square root of the
modular function (see e.g. \cite{cowling-giulini-etal}), by means of
refinements of the estimates for spherical functions in \cite{ionescu}.
However, we shall not pursue this here, since we prefer to present the
completely self-contained proof which we can give for the affine group.

In Section \ref{ints}, we prove a lemma which describes integration of
radial functions
over the affine groups. In Section \ref{exps}, we derive an explicit
kernel representation for the resolvent of $L$ using the theory of
hypergeometric functions. It is known (\cite{hulanicki}, \cite{gnewuch}) that
the resolvent kernels can be expressed in terms of special functions, and
this has been used in \cite{gqs} to prove estimates for singular integrals
related to $L$. We chose
to derive the integral formula for the resolvent kernel from scratch,
even though this task could have been done quoting tables of special
functions such as \cite{magnus}. We hope some readers will find benefit of
our explicit calculations.

Section 4 presents a subordination argument to obtain convenient integral
representations for the convolution kernels of the operators (\ref{epsi}).

In Section \ref{refine}, we prove some asymptotic formulas for
the hypergeometric functions appearing in Section \ref{exps}.

Section \ref{specloc} assembles the results of the previous sections
to prove Theorem \ref{ktsatz}, which states pointwise estimates for
the convolution kernel of (\ref{epsi}), and Proposition \ref{wlambdat},
which states $L^1(G)$ estimates for these kernels. We
also reprove a multiplier theorem of \cite{hebisch-steger}.

In Section \ref{sobolev}, we prove growth estimates for the wave
propagator associated with $L$ using spectrally defined Sobolev norms.

The second author would like to thank Adam Sikora and Terry Tao for
helpful discussions about the affine group and wave operators in general.

\section{Integration of radial functions on the affine group} \label{ints}

In this section we discuss integration of ``radial'' functions. The results
 will be useful in the  estimation of $L^1$-norms of convolution kernels
for functions of the Laplacian $L$.

Bending the notion of radial function, by radial function we mean a
function of the type
$$e^{-n x/2}g(R)$$
with $R$ as in (\ref{defr}).

We briefly motivate the special form of the radial variable $R$. For $n=1$,
the affine group $G$ is a subgroup of the group of conformal automorphisms of
the upper half plane via the identification of $(x,y)$ with the map $z\to
e^x z + y$. This subgroup acts simply transitively on the upper half plane,
thus we can naturally identify $G$ as a set with the upper half plane,
identifying the neutral element with the point $i$.
In particular, the hyperbolic metric on the upper half plane turns out
left-invariant. The pull back of the hyperbolic distance from a point $z$ to
the point $i$ (which is $\log|(1+\rho)/(1-\rho)|$ with $\rho=|z-i|/|z+i|$)
gives a natural ``radial'' distance from the origin in the affine group given
by
$$\arccosh (\cosh x+ \frac 12 y^2 e^{-x}).$$
Since $R(x,y)=R((x,y)^{-1})$
there is no difference between a left and a right radial variable.

\begin{lemma}\label{intlem}
Given a function $g:\R^+\to \C$, then
$$\int_G e^{-n x/2 } g(R(x,y))\, dxdy= \int_0^\infty g(R) J(R)\, dR,$$
where
\begin{eqnarray}
J(R)\sim R^n, &&\text{ if} \  R\le 1,  \\
J(R)\sim R e^{nR/2},&&\text{ if} \  R\ge 1,
\end{eqnarray}
and $a\sim b$ means that each of the two numbers can be bounded by a
constant times the other, the constant depending only on $n$.
\end{lemma}

\begin{proof}

Define
$$B(r)=\int_{R(x,y)\le r} e^{-nx/2}\, dx dy.$$ Then we have for all $r\ge 0$
$$J(r)=B'(r).$$
Thus we have to estimate $B'(r)$.
Observe that $R(x,y)\le r$ implies $x\le r$ and $$\|y\|\le
(2e^x (\cosh (r)-\cosh(x)))^{1/2}.$$
Doing the $y$-integration and letting $V_n$
be the Euclidean volume of the unit ball in $\R^n$, we obtain $$B(r)=V_n
\int_{-r}^r e^{-nx/2} (2e^{x} (\cosh (r)-\cosh(x)))^{n/2}\, dx.$$
Simplification and differentiation gives $$B'(r)=\frac n2 2^{n/2} V_n
\int_{-r}^r
(\cosh (r)-\cosh(x))^{n/2-1}\sinh(r) \, dx. $$

First assume $r\ge 1$. We break the integral in the previous display
into the sum of
\begin{equation}\label{firsti}
I_1= \sinh(r)\int_{|x|<r-1}
(\cosh (r)-\cosh(x))^{n/2-1}\, dx
\end{equation}
and
\begin{equation}
\label{secondi}
I_2=\sinh(r)
\int_{r-1\le |x|\le r}
(\cosh (r)-\cosh(x))^{n/2-1}\, dx\ .
\end{equation}
In the first integral, we have $|x|<r-1$ and thus $$\cosh R
-\cosh(x)\sim e^r\ ,\ \ \ \sinh(r)\sim e^r.$$ Therefore
$$I_1\sim r e^{nr/2} .$$
Thus $I_1$ has already the correct order of magnitude which we need to show
for $I_1+I_2$. Since $I_1$ and $I_2$ are positive, we only need an upper bound
for $I_2$. Since in the domain of integration of
$I_2$ we still have
$$\cosh (r)-\cosh(x)\le 2 e^r\ ,\ \ \ \sinh(r)\le e^r\ ,$$ we can do the same
calculation as before to obtain $$I_2\le C e^{rn/2}.$$

Now assume $r<1$. We do a similar splitting of the integral as before, now
into the regions $|x|<r/2$ and $r/2\le |x|\le r$. Call the corresponding
integrals $I_1$ and $I_2$. In the domain $|x|<r/2$ we have
$$\cosh(r)-\cosh(x) \sim r^2\ , \ \ \ \sinh(r)\sim r\ , \ \ \ e^x\sim 1.$$
Hence
$$I_1 \sim r^{n-1} \int_0^{r/2}\, dx\sim r^n.$$
Similarly as before, in the
domain $r/2\le |x|\le r$ we have the same upper bounds as in the domain
$|x|<r/2$, and thus obtain $I_2\le C r^{n}$. This completes the proof of Lemma
\ref{intlem}.
\end{proof}

\section{An explicit kernel for the resolvent of $L$} \label{exps}

Assume that $k$ is an integrable function on $G$ such that for
every compactly supported smooth function $\varphi$ we have (in the
distributional sense, for this purpose we identify $G$ with $\R^{n+1}$)
\begin{equation}\label{fundsol}
\int \varphi(g)(L-\lambda)
k(g)\, dg=\varphi(0)\ .
\end{equation}
Then, by left invariance of $L$ and right invariance of
$dg$, $$(L-\lambda) \int \varphi(g^{-1})k(gg')\, dg= \int
\varphi(g^{-1})[(L-\lambda) k](gg')\, dg$$ $$=
\int \varphi(g'g^{-1})[(L-\lambda) k](g)\, dg= \varphi(g').$$
Thus the resolvent operator $(L-\lambda)^{-1}$ is given by right
convolution with $k$, which extends to a bounded operator on $L^2(G)$. The
following lemma describes such a fundamental solution $k$.

\begin{lemma}\label{resolvent_kernel_lemma}
 Assume $\lambda\in \C\setminus [0,\infty)$ and choose
$\nu:=i\sqrt{\lambda}$ with strictly negative real part.
Then the resolvent operator $(L-\lambda)^{-1}$ is given by right
convolution with the kernel $k$ defined by ($R=R(x,y)$ as in (\ref{defr})):

$$k(x,y)= (-1)^l \frac{2^{-1-n/2}\pi^{-n/2}}
{\Gamma(1-n/2+l)}e^{-nx/2}
\int_{R}^\infty D_{\sinh,v}^l[e^{\nu v}](\cosh v-\cosh R)^{-n/2+l} \, dv\ ,$$
where $l$ is any integer with $-n/2+l>-1$
and we have written $D^l_{\sinh,v}$ for the $l$- th power of
$D_\sinh: g\to D(g/\sinh)$
acting in the $v$ variable.

Moreover, the kernel $k$ satisfies the estimate
\begin{equation}\label{locint}
 \int_{B_R} |k| \, dxdy \le  C_n (1+|\nu|)^{\frac n2}
[1 + \int_0^R e^{\Re(\nu) r} r\, dr]\ ,
\end{equation}
where $B_R$ is the ball of radius $R$ about the origin in $G$.
In particular, $k\in L^1(G)$.

\end{lemma}

{\it Proof.}
Fix $\lambda$.
For $n>1$, we will show that $k$ as defined in the lemma is integrable and
smooth outside
the origin, satisfies $(L-\la) k=0$ outside the origin, and has asymptotics
$$k(x,y)=2^{-2}\pi^{-(n+1)/2} \Gamma((n-1)/2) (x^2+\|y\|^2)^{-(n-1)/2}+
O(x^2+\|y\|^2)^{- n/2}$$
near the origin. This will give for any compactly supported $\varphi$
$$\int \varphi(g) (L-\la)k(g)\, dg = \int \eta(g)\varphi(g) (L-\la) k(g)\,
dg\ ,$$
where $\eta$ is a smooth cutoff  function at scale $\ve,$ i.e. constant equal
to $1$ on an $\ve$- neighborhood around the
origin and  zero  outside a $2\ve$ neighborhood, with the usual contol of
derivatives. By definition of the distributional derivative, the last
display becomes $$\int (L-\la) (\eta \varphi)(g) k(g)\, dg\ .$$
If we subtract the leading order term from $k$ in this integral, then the
remaining integral tends to zero as $\ve\to 0$. Thus we may replace $k$ by
the leading order term. Also, we may disregard in the expansion of
$(L-\la)(\eta \psi)$ by Leibniz' rule all terms other than those taking two
derivatives of $\eta$, and also may approximate the coefficient $e^x$ by
$1$. Thus the last display is equal to ($\Delta=- \sum\partial_j^2$) $$
2^{-2}\pi^{-(n+1)/2} \Gamma((n-1)/2)
\int \Delta (\eta\varphi) \left(x^2+|y|^2 \right)^{-(n-1)/2}\, dxdy\ .$$ Now
standard theory in $R^m$ (\cite{taylor} pp. 211, 262) gives that the
last display is equal to $\varphi(0)$, which was to be proved. If $n=1$, we
use the same approach, here the asymptotic behaviour of $k$
near the origin is
$$2^{-2}\pi^{-1} \log(x^2+\|y\|^2)+O(1)\ .$$

Now we prove the properties of $k$ claimed above. Define
$$d(x,y):=\cosh(R(x,y))=\cosh(x)+\frac 12 \|y\|^2 e^{-x}\ .$$ Then the kernel
$k$ is of the form
\begin{equation}
\label{definek}
k(x,y)=e^{-nx/2}f(d(x,y)),
\end{equation}
with a function $f$ which is smooth on $(1,\infty)$. We claim that
$(L-\lambda)k=0$ outside the origin is equivalent to $f$ satisfying the
ordinary differential equation \begin{equation}\label{ode}
-\frac {n^2}4 f(d)- (n+1) d f'(d)- (d^2-1) f''(d)=\lambda f(d) \end{equation}
for $d>1$.
To verify the claim, we observe
$$(Xd)^2+\sum_{j=1}^n (Y_j d)^2$$
$$=(\frac 1 2 e^x -\frac 12 e^{-x}-\frac 12 \|y\|^2e^{-x})^2+\|y\|^2 = (\frac
1 2 e^x +\frac 12 e^{-x}+\frac 12 \|y\|^2e^{-x})^2-1$$ $$=d^2-1$$
and
$$ X^2 d +\sum_{j=1}^n Y_j^2d - n X d$$
$$=(\frac 12 e^x +\frac 12 e^{-x}+ \frac 12 \|y\|^2 e^{-x})+ ne^x -n (
\frac 12 e^x - \frac 12 e^{-x}- \frac 12 \|y\|^2 e^{-x}) =(n+1)(\frac 12
e^x +\frac 12 e^{-x}+ \frac 12 \|y\|^2 e^{-x})$$ $$= (n+1) d\ .$$
Therefore we have
$$(X^2+\sum_{j=1}^n Y_j^2) e^{-nx/2} f(d(x,y))=$$ $$=e^{-nx/2}
[((Xd)^2+(Yd)^2)f''(d)+(X^2d+Y^2d)f'(d)-n(Xd)f'(d)+\frac {n^2}4 f(d)]$$
$$=e^{-nx/2}[\frac{n^2}4 f(d(x,y))+
(n+1) d(x,y) f'(d(x,y))+ (d(x,y)^2-1)f''(d(x,y)]\ .$$ Thus if $f$ satisfies
the ordinary differential equation (\ref{ode}) on $(1,\infty)$, then
$(L-\lambda)k=0$ outside the origin, and conversely. Equation (\ref{ode})
is a classical
hypergeometric differential equation. The Riemann symbol \cite{klein}
associated to this differential equation is
\begin{equation}
P\left(\begin{array}{cccc} -1 & 1 & \infty &   \\
0 & 0 & \frac n2 + \nu & d \\
-\frac {n-1}{2} & -\frac{n-1}2 & \frac n2 - \nu &
\end{array} \right)\ .
\end{equation}
There is a two dimensional space of solutions $f$. However, there is only
a one dimensional space of solutions (those with leading asymptotics
$d^{-n/2+\nu}$ as $d\to \infty$), which make $k$ as defined above
integrable on the group $G$.
Of course exactly one of these solutions is normalized properly to make $k$ a
fundamental solution.

The following lemma provides an explicit solution of the
differential equation $(\ref{ode})$ in a certain complex
region of parameters $n$ and $\nu$.
\begin{lemma}\label{fzerolemma}
Assume $-\Re(n)/2>-1$ and $\Re(\nu)-\Re(n)/2<0$. Then the function
\begin{equation}\label{fzero}f_0(d)=
\int_{\arccosh (d)}^\infty e^{\nu v}(\cosh v-d)^{-n/2} \, dv\ ,
\end{equation}
defined for $d>1$,
satisfies the ordinary differential equation (\ref{ode}). \end{lemma}

\begin{proof}
Under the stated assumptions on $n$ and $\nu$, the integral defining $f_0$
is absolutely integrable. We first assume that $-\Re(n)/2>1.$

 Differentiating under the integral
sign gives $${f_0}'(d)= \frac n2 \int_{\arccosh(d)}^\infty
e^{\nu v} (\cosh v - d)^{-n/2-1}\, dv\ ,$$
$${f_0}''(d)= \frac n2(\frac n2 + 1) \int_{\arccosh(d)}^\infty
e^{\nu v} (\cosh v - d)^{-n/2-2}\, dv\ .$$ Hence,
$$(n+1)d{f_0}'= (\frac{n^2} 2+\frac n2)(- f_0 + \int_{\arccosh(d)}^\infty
e^{\nu v} (\cosh v - d)^{-n/2-2} (\cosh v^2-d\cosh v)\, dv\ ,$$
$$(d^2-1){f_0}''= (\frac {n^2}4 +\frac n 2) (f_0 -
\int_{\arccosh(d)}^\infty
e^{\nu v} (\cosh v - d)^{-n/2-2} (\cosh v^2-2d\cosh v+1)\, dv\ .$$ It thus
remains to prove
$$
\int_{\arccosh(d)}^\infty
e^{\nu v} (\cosh v - d)^{-n/2-2} \left[-\frac{n^2}4 \cosh v^2 - \frac n2
d\cosh v + \frac {n^2} 4+\frac n2 \right]\, dv$$ $$=\lambda
\int_{\arccosh(d)}^\infty
e^{\nu v} (\cosh v - d)^{-n/2} \, dv\ .$$
However, by partial integrations, the right hand side is equal to
$$=\frac{\lambda n}{2\nu}
\int_{\arccosh(d)}^\infty
e^{\nu v} (\cosh v - d)^{-n/2-1} \sinh(v) \, dv$$ $$=-\frac{\lambda n}{2\nu^2}
\int_{\arccosh(d)}^\infty
e^{\nu v} (\cosh v - d)^{-n/2-1} \cosh v \, dv$$ $$+\frac{\lambda}{\nu^2}
(\frac{n^2}4 +\frac n2) \int_{\arccosh(d)}^\infty
e^{\nu v} (\cosh v - d)^{-n/2-2} \sinh(v)^2 \, dv\ .$$ This proves the
lemma for
$-\Re(n)/2>1$
since $\lambda/\nu^2=-1$ and $\sinh^2(v)=\cosh^2(v)-1$

The case $-\Re(n)/2>-1$ then follows by analytic continuation.
\end{proof}

Clearly, any analytic continuation of (\ref{fzero}) in the parameters $n$
and $\nu$
also satisfies the ordinary differential equation (\ref{ode}). The following
lemma provides explicit expressions for such analytic continuations.
\begin{lemma}\label{fzerol}
Assume $-\Re(n)/2>-1$ and $\Re(\nu)-\Re(n)/2<0$. Then, for each integer $l\ge
-1$, we have the identity $$\Gamma(- \frac n2+1)^{-1}f_0(d)= (-1)^l
\Gamma(l-\frac n2+1)^{-1} \int_{\arccosh d}^\infty D^l_{\sinh,v} [e^{\nu
v}] (\cosh v-d)^{-n/2+l} \, dv\ .$$
The right hand side provides
the unique analytic continuation of the left hand side to the parameter region
$-\Re(n/2)+l>-1$ and $\Re(\nu)-\Re(n/2)<0$. The right hand side satisfies the
ordinary differential equation (\ref{ode}) in this parameter region.
Here we have set for $\Re(\nu)<0$
$$D^{-1}_{\sinh,v}[e^{\nu v}]= \nu^{-1}\sinh(v) e^{\nu v}\ .$$

\end{lemma}

\begin{proof}
By partial integration, we have for $l\ge 0$ $$
\Gamma(l-n/2+1)^{-1}
\int_{\arccosh d}^\infty D^l_{\sinh,v} [e^{\nu v}] (\cosh v-d)^{-n/2+l} \, dv$$
$$
= - \Gamma(l-n/2+1)^{-1}
\int_{\arccosh d}^\infty D^{l-1}_{\sinh,v} [e^{\nu v}] \sinh(v)^{-1}D_v[
(\cosh v-d)^{-n/2+l}] \, dv$$ $$
= - \Gamma(l-n/2)^{-1}
\int_{\arccosh d}^\infty D^{l-1}_{\sinh,v} [e^{\nu v}] (\cosh
v-d)^{-n/2+l-1} \, dv$$
By induction, this proves the identity of the lemma.
\end{proof}

To normalize the resolvent kernel properly, we need to calculate the
asymptotic behaviour of $f_0$ near $1$. The reader not interested in
explicit constants may skip the following lemma.

\begin{lemma}\label{fconstants}
Let $\Re(\nu)-\Re(n/2)<0$ and assume $n$ is not
an odd negative integer. For $R>0$ and $R$ near $0$ and any $\ve\in (0,1)$
we have
$$\Gamma(1-n/2)^{-1} f_0(\cosh R)= 2^{n/2-1}\pi^{-1/2} \Gamma((n-1)/2)
R^{1-n}
+O(R^{1-\Re(n)+\ve})$$
if $n\neq 1$ and
$$
\Gamma(1-n/2)^{-1}f_0(\cosh R)= 2^{1/2}\pi^{-1/2} |\log(R)| + O(1)
$$
in case $n=1$. Here the left hand side is to be understood as an analytic
function in the sense of Lemma \ref{fzerol}. \end{lemma}

\begin{proof}
We do the case $1<\Re(n)<2.$ The general case follows by
similar calculations or by methods of analytic continuation.

Assume $R<<1$.
We split the integral
$$f_0(\cosh R)= \int_{R}^\infty e^{\nu v}(\cosh v-\cosh R)^{-n/2} \, dv$$
into the integrals over the intervals $[R,R^{1-\ve}]$ and
$[R^{1-\ve}, \infty]$.
Since $\cosh v-1> v^2$, the integral over the second interval is bounded by
$$C \int _{R^{1-\ve}}^\infty v^{-\Re(n)} \, dv \le C
R^{(1-\Re(n))(1-\ve)}\ .$$ Thus this integral is negligible. For
$v<R^{1-\ve}$ we write $$e^{\nu v}=1 + O(R^{1-\ve})\ ,$$
$$\cosh v-\cosh R= \frac 12 (v^2-R^2)(1+O(R^{2-\ve}))\ .$$ Thus the first
integral is
$$ 2^{n/2} \int_{R}^{R^{1-\ve}} [v^2-R^2]^{-n/2} \, dv
(1+O(R^{1-\ve}))\ .$$
By an argument as before, we can change the domain of integration back to
$[R,\infty]$ doing at most an error of order
 $R^{(1-n)(1-\ve)}$. Thus we have to get asymptotics of the integral $$
2^{n/2}\int_{R}^{\infty} [v^2-R^2]^{-n/2} \, dw$$ $$=2^{n/2-1}
R^{1-n}\int_{1}^{\infty} [w-1]^{-n/2} w^{-1/2} \, dw\ .$$ The integral in
the last display can be expressed in terms of Gamma functions:
$$\int_{1}^{\infty}
[w-1]^{-n/2}w^{-1/2} \, dw$$ $$= \Gamma(1/2)^{-1}\int_{1}^{\infty}
[w-1]^{-n/2} \int_0^\infty t^{1/2} e^{-tw} \, \frac{dt}t \, dw$$ $$=
\pi^{-1/2}\int_0^\infty t^{1/2} e^{-t} \int_{0}^{\infty} r^{-n/2} e^{-tr} \,
dr \, \frac{dt}t $$ $$= \pi^{-1/2}\Gamma(1-n/2) \int_0^\infty t^{n/2-1/2}
e^{-t} \, \frac{dt}t $$
$$= \pi^{-1/2} \Gamma(1-n/2)\Gamma((n-1)/2)\ .$$
This proves the asymptotics claimed in the lemma.
\end{proof}

With Lemma \ref{fconstants} we have completed the proof of the identities
for the resolvent kernel $k$ claimed in Lemma
\ref{resolvent_kernel_lemma}.
\medskip

It remains to prove the $L^1$ estimates
for $k$ (compare with \cite{hulanicki}). We can assume that $n$ is a positive
integer.
We have to estimate the integral
$$\int_R^\infty D_{\sinh,v}^l[e^{\nu v}] (\cosh(v)-\cosh(R))^{l-n/2}\, dv
=H^1_\nu(R)+H^2_\nu(R),$$ where
$$H^1_\nu(R)=\int_R^{R+2} D_{\sinh,v}^l[e^{\nu v}]
 (\cosh(v)-\cosh(R))^{l-n/2}\, dv$$
and
$$H^2_\nu(R)=\int_{R+2}^\infty D_{\sinh,v}^l[e^{\nu v}]
 (\cosh(v)-\cosh(R))^{l-n/2}\, dv,$$ and
where we choose $l$ such that $l-n/2\in \{-1/2,0\}.$

It is easily seen that
$$|H^2_\nu(R)|\lesssim |\nu|^l\int_{R+2}^\infty e^{(-l+\Re(\nu))v}\,
e^{(l-\frac n2)v}\, dv\lesssim(1+|\nu|)^{\frac n2}e^{-\frac n2 R+\Re(\nu)R}.$$

To estimate $H^1_\nu(R),$ for $R\ge 1$ we use that in the domain of
integration
we have $$D^m \sinh^{-1}(v)\sim e^{-R},$$
$$|D^m e^{\nu v}|\sim |\nu|^m e^{\Re(\nu) R},$$
$$\cosh(v)-\cosh(R)\sim e^{R}(v-R).$$

This leads to the same estimate as for $H^2_\nu(R).$

For $R<1$ we use that in the domain of integration
 $$D^m \sinh^{-1}(v)\sim v^{-1-m},$$
$$\cosh(v)-\cosh(R)\sim v(v-R),$$
$$|D^m e^{\nu v}|\sim |\nu|^m .$$
Thus $$|H^1_\nu(R)|\lesssim (1+|\nu|^l) \int_R^{R+2}
v^{-l-n/2}(v-R)^{-n/2+l}\, dv$$
 $$\lesssim (1+|\nu|^l) R^{1-n} \int_1^{\infty} v^{-l-n/2}(v-1)^{-n/2+l}\,
dv$$
$$\lesssim (1+|\nu|^{\frac n2}) R^{1-n} $$

Put together, we find that for  $\Re(\nu)<0$
$$\left |\int_R^\infty D_{\sinh,v}^l[e^{\nu v}] (\cosh v-\cosh R)^{l-n/2}\, dv
\right |
\le C_n (1+|\nu|)^{\frac n2} R^{1-n} (1+R^{n-1}) e^{-\frac n2 R + \Re(\nu)
R}\ .
$$
Using Lemma \ref{intlem}, this proves the desired estimate (\ref{locint})
and completes the proof of Lemma \ref{resolvent_kernel_lemma}.

\section{Spectral Multipliers}\label{specmul}
Let $\psi\in C_0({\mathbb R}).$ Since $\psi(L)$ is a bounded linear operator
on $L^2(G),$ by the Schwartz' kernel theorem and left-invariance of $L,$ there
exists
a unique convolution kernel  $k_\psi\in\mathcal D'(G)$ such that
$\psi (L)\vphi=\vphi\ast k_\psi,\ \vphi\in\mathcal D(G).$ We shall derive
an integral
representation for $k_\psi$ in this section. In the sequel, we shall
sometimes also use
the suggestive notation  $k_\psi=\psi(L)\delta_0.$

 Since the
Gauss-kernels
$g_{\varepsilon}(s)=\frac{1}{\sqrt{2\pi\ve}}e^{-\frac{s^2}{2\ve}}$,
$\ve
>0$, form an approximation to the identity with respect to convolution, we have
$\psi_{\ve}:=\psi\star g_{\ve}\rightarrow \psi$, uniformly as
$\ve\rightarrow 0$. But, $\psi_{\ve}$ has an analytic continuation, given
by
$$\psi_{\ve}(\zeta):=\frac{1}{\sqrt{2\pi\ve}}\int_{\mathbb R}\psi
(t)e^{-\frac{(t-\zeta)^2}{2\ve}}dt,\ \zeta\in{\mathbb C}\, .$$ Therefore,
$\psi_{\ve}(L)$ is given by the Cauchy-integral
$$\psi_{\ve}(L)=\frac{1}{2\pi
i}\int_{\gam_\delta}\psi_{\ve}(\zeta)(L-\zeta )^{-1}\,d\zeta\, ,$$
where, for any $\delta>0$, we may choose for $\gam_{\delta}$ the
(clockwise) contour
$\gam_{\delta}:s\mapsto (s+i\delta)^2,$ $s\in{\mathbb R}$. By
Lemma \ref{resolvent_kernel_lemma},
$\psi_{\ve}(L)\vphi =\vphi\star k{_{\psi_{\ve}}}$,
$\vphi\in\mathcal{S}$, where the kernel $k{_{\psi_{\ve}}}$ is given, for any
$l>\frac{n}{2}-1$, by
\begin{equation}\label{kpeps}
k{_{\psi_{\ve}}}(x,y)=c_le^{-\frac{nx}{2}}\int_{\mathbb
R}\psi_{\ve}((s+i\delta)^2)F_R(s+i\delta)(s+i\delta)\, ds\, ,
\end{equation}
where
\begin{equation}\label{cl}
c_l:=\frac{(-1)^l}{\pi i}\frac{2^{-1-\frac{n}{2}}
\pi^{-\frac{n}{2}}}{\Gamma (-\frac{n}{2}+1+l)}\, ,
\end{equation}
and where
\begin{equation}\label{FRzeta}
F_R(\zeta) :=\int^{\infty}_R D^l_{\sinh ,v}[e^{i\zeta v}](\cosh v-\cosh
R)^{-\fu2 +l}\,dv\, ,
\end{equation}
if $\Im\zeta >-\frac{n}{2}$ .\\

The estimate (\ref{locint})  shows that the mappings $(x,y)\mapsto
e^{-\frac{nx}{2}}F_{R(x,y)}(s+i\delta)$ are locally integrable
on $G$, and their integrals over compact subsets of $G$ are of
polynomial growth in $s$, uniformly in $0<\delta<1$.

Moreover, since
$$\Re[(t-(s+i\delta)^2)^2)] - s^4 \le C [1+|s|]^3\ ,$$
where $C$ is uniform in $0\le \delta\le 1$, $t\in \supp(\psi)$,
we have
$$\mid\psi_{\ve}((s+i\delta)^2)\mid\le C
e^{-cs^4}\ ,$$
where $C$ and $c>0$ depend on $\ve$ but not on $\delta$.

Therefore, given $\vphi\in
\mathcal{D}(G)$, by the dominated convergence theorem the limit of
$$\int \psi_{\ve}((s+i\delta)^2)
{\dst\int_G}e^{-\frac{nx}{2}}F_{R(x,y)}(s+i\delta)\,
\vphi(x,y)dxdy\,(s+i\delta) ds$$
as $\delta$ tends to $0$ is equal to the same expression with $\delta=0$.
Therefore, in the sense of distributions,
\begin{equation}\label{kpsieps}
k{_{\psi_{\ve}}}(x,y) = c_le^{-\frac{nx}{2}}\int_{\mathbb R}
\psi_{\ve}(s^2) F_{R(x,y)}(s)s\,ds\, .
\end{equation}

Finally, as $\ve \rightarrow 0$, $\psi_{\ve}(L) \rightarrow \psi(L)$ in the
operator norm on $L^2(G)$, which implies that
$k_{\psi_{\ve}} \rightarrow k_{\psi}$ in $\mathcal{D}'(G)$. On the other hand,
$|\psi_{\ve}(s^2)| \leq C e^{-\frac{s^4}{4}}$,  for $0< \ve < 1$, so that,
again by the dominated convergence theorem and (\ref{kpsieps}),
 $k_{\psi_{\ve}} \rightarrow C_{l} e^{-\frac{nx}{2}} \int_{\mathbb R} \psi(s^2)
F_R (s) s ~ds$, in the sense of distributions. We have thus proved
\begin{proposition}\label{nullunend}
Let $\psi \in C_0(\mathbb R).$ Then,
for any $l>\frac{n}{2}-1$, the convolution kernel $k_{\psi}$ of $\psi(L)$
is locally integrable on $G$, and is given by
\begin{eqnarray}\label{sub}
k_{\psi}(x,y) =&& c_l e^{-\frac{nx}{2}}
\int_{\mathbb R} \psi(s^2) F_{R(x,y)}(s) \, s \,ds\nonumber\\
=&&
 c_l e^{-\frac{nx}{2}}
\int_{0}^\infty \psi(s^2)[F_{R(x,y)} (s)-F_{R(x,y)} (-s)] s \,ds,
\end{eqnarray}
with $c_l$ given by
(\ref{cl}).\end{proposition}

\section{Asymptotics  of $F_R(s)$}\label{refine}

We denote by $S^{\alpha}$ the  symbol class
$$  S^{\alpha}:=\{b\in C^{\infty}(\mathbb R) :\parallel b
\parallel_{S^{\alpha},k}:=
\sup_s (1+s^2)^{\frac{-\alpha +k}{2}}\mid b^{(k)}
(s)|~<\infty\,\ \text{for all}\ k\in\mathbb N\} . $$
The spaces $S^\alpha$  are Fr\'echet-spaces, with the topology
induced by the sequence of semi-norms $\parallel
\cdot\parallel_{S^{\alpha},k},\, k\in\mathbb N$.
The product of a function in $S^\alpha$ with a function in $S^{\beta}$
is in $S^{\alpha+\beta}$. Moreover, $S^\alpha\subset S^\beta$
if $\alpha<\beta$ and $Db\in S^{\alpha-1}$ if $b\in S^\alpha$.
The following general lemma will also be useful:
\begin{lemma}\label{movers}
Let $b_\alpha\in S^\alpha$ and $\beta,\gamma>0$.
Then we can write for each $R>0$
$$b_\alpha = R^\beta b_{\alpha+\beta,R} + R^{-\gamma} b_{\alpha-\gamma,R}$$
with $b_{\alpha+\beta}\in S^{\alpha+\beta}$ and
$b_{\alpha-\gamma}\in S^{\alpha-\gamma}$ uniformly in $R$, i.e.,
$$\|b_{\alpha+\beta}\|_{S^{\alpha+\beta},k} \le C_k$$
$$\|b_{\alpha-\gamma}\|_{S^{\alpha-\gamma},k} \le C_k$$
with constants $C_k$ independent of $R$.
\end{lemma}
Note: We have suppressed the $R$-dependence of the
symbols $b_{\alpha+\beta}$, $b_{\alpha-\gamma}$ in the notation,
and we will continue to suppress any $R$ dependence of symbols $b$
throughout the rest of this paper.
\begin{proof}
Let $\chi$ be a smooth cutoff function which is constant $1$
on $(-\infty, 1]$ and vanishes on $[2,\infty)$.
Then we write
$$b_\alpha= R^\beta\,[R^{-\beta}(1-\chi(R(1+s^2)^{1/2}))b_\alpha]+
R^{-\gamma}\,[R^\gamma \chi(R(1+s^2)^{1/2})b_\alpha] \ .$$
It is easy to see by Leibniz' rule that this is the desired
splitting.
\end{proof}

We wish to estimate the function
$$ F_R(s):=\int\limits^{\infty}_{R}
D^l_{\sinh, v}[e^{isv}](\cosh v-\cosh R)^{-\frac{n}{2}+l}dv,\quad s\in\R,
\quad(l>\frac{n}{2}-1). $$
The estimates are stated in Proposition \ref{r1} for the case $R\ge 1$ and
in Proposition \ref{o<Rle1} for the case $0<R\le 1$.

\begin{proposition}\label{r1}
If $R\ge 1$, then
\begin{equation}\label{frs}
F_R(s)=e^{-\frac{n}{2}R}\, e^{iRs}\, b_{\frac n2-1}(s),
\end{equation}
where $b_{\frac n2-1}\in S^{\frac n2-1}$ uniformly in $R$.
\end{proposition}

This proposition will follow from the subsequent lemmas.
\begin{lemma}\label{dshiv}
For $v\ge 1$ we can write
$$D^l_{\sinh ,v}[e^{isv}]={\dst\sum^l_{k=0}} s^k \, q_k(v)\, e^{-lv}\,
e^{isv}$$
with $ q_k \in S^0$ for all $k$.
\end{lemma}

\begin{proof} This is proved by induction on $l\in\mathbb N$, the case
$l=0$ being
trivial. Assume the statement is true for some $l\in\mathbb N$. Then
$$ D^{l+1}_{\sinh,v} [e^{isv}]= \sum^l_{k=0}\, s^k\, D_v
( (\sinh v)^{-1} \, {q_k(v)}\, e^{-lv} \,e^{isv}) $$
$$  ={\dst\sum^l_{k=0}} \,s^k \, D_v( \frac{2}{1-e^{-2v}}\,
{q_k(v)}\,
e^{-(l+1)v}\, e^{isv}). $$
On the interval $[1,\infty)$, the function
$\frac{1}{1-e^{-2v}}={\dst\sum^{\infty}_{m=0}}
 e^{-2mv}$ coincides with
a function in $S^0$. This easily implies the statement of the lemma for $l+1$.
\end{proof}

We can therefore decompose

\begin{eqnarray*} F_{R}(s)&=&
\sum^l_{k=0} \,
s^k\int^{\infty}\limits_{R}\;
q_k(v)(\cosh v-\cosh R)^{-\frac{n}{2}+l} e^{-lv} e^{isv}\, dv \\
&=&
\sum^l_{k=0} \,
e^{iRs}e^{-\frac{n}{2} R}s^k\int^{\infty}\limits_{0}
q_k(R+v)\left[(\cosh
(R+v)-\cosh R)e^{-(R+v)}\right]^{-\frac{n}{2}+l} e^{-\frac n2 v}
e^{isv}\, dv \ .\\\end{eqnarray*}
Fixing $k\in\{0,\ldots ,l\}$ and writing
$$ \gam_R(v):=
q_k(R+v)\left[(\cosh
(R+v)-\cosh R)e^{-(R+v)}\right]^{-\frac{n}{2}+l} e^{-\frac n2 v}\ ,$$
it then suffices to
prove that the function
\begin{equation}\label{frsint}
f_R(s):=\int^{\infty}\limits_{0}\gam_R(v)e^{isv} dv,\quad s\in\mathbb R ,
\end{equation}
lies in $S^{\frac{n}{2} -l-1}$ uniformly with respect to  $R\ge 1$.\\

Let $\chi$ be a smooth cutoff function which is constant
equal to $1$
on $(-\infty,1]$ and vanishes on $[2,\infty)$.

It suffices to show that
\begin{equation}\label{frsint1}
f_{R,1}(s):=\int^{\infty}\limits_{0} \chi(v) \gam_R(v)e^{isv} dv,\quad
s\in\mathbb R ,
\end{equation}
\begin{equation}\label{frsint2}
f_{R,2}(s):=\int^{\infty}\limits_{0} (1-\chi(v)) \gam_R(v)e^{isv} dv,\quad
s\in\mathbb R ,
\end{equation}
are in $S^{\frac{n}{2} -l-1}$ uniformly with respect to  $R\ge 1$.\\

The following lemma settles the question for $f_{R,2}$.

\begin{lemma}\label{globg}
The function
$(1-\chi) \gam_R$ is in the Schwartz class uniformly in $R>0$.
\end{lemma}
\begin{proof}
Since the function $q_k$ is in $S^0$, all derivatives $D_v^l q_k$
are bounded.
It then suffices to show that also all derivatives of
\begin{equation}\label{coshpow}
\left[(\cosh
(R+v)-\cosh R)e^{-(R+v)}\right]^{-\frac{n}{2}+l}
\end{equation}
are bounded on $[1,\infty)$, uniformly in $R>0$.
However,
\begin{equation}\label{addth}
(\cosh
(R+v)-\cosh R)e^{-(R+v)}
=\frac{[\cosh(v)-1]}{e^v} \frac {\cosh R}{e^R}
+\frac{\sinh(v)}{e^v} \frac {\sinh R}{e^R},
\end{equation}
and this function and all its derivatives are
bounded on $[1,\infty)$ uniformly in $R>0$.
Moreover, (\ref{addth}) is also bounded below by some
$\ve>0$ on $[1,\infty)$ uniformly in $R>0$.
Therefore, all derivatives of (\ref{coshpow})
are bounded, which completes the proof of the lemma.
\end{proof}

It remains to show that $f_{R,1}$ is in $S^{\frac n2 -l-1}$.
This will follow from the next two lemmas.

\begin{lemma}\label{locg}
For $v>0$ we can write
$$\chi \gam_R = g_R(v) v^{-\frac n2 + l}\ ,$$
where $g_R$ is supported in $v\le 2$ and
$D^kg_R$ is bounded uniformly in $R\ge 1$
for all $k\in \N$.
\end{lemma}

\begin{proof}
Taylor expansion of $\sinh(v)$ and $\cosh(v)$
in the expression (\ref{addth}) gives
for $v\le 2$
$$(\cosh
(R+v)-\cosh R)e^{-(R+v)}= v \tilde{g}(v)$$
for some function $\tilde{g}$ which is bounded below
by $\ve>0$ and has all derivatives bounded above
uniformly in $R\ge 1$.
This proves the lemma.
\end{proof}

\begin{lemma}\label{gam}
Let $g\in S^0$ be supported in $[-2,2]$ and $\alpha>-1.$
Then
$$f (s):= \int_0^\infty g(v) v^\alpha e^{ivs}\, dv, \quad s\in\R,$$
is in $S^{-\alpha-1}$
with semi-norms $\parallel f\parallel_{S^{-\alpha -1},k}$ controlled by the
seminorms $\parallel g \parallel_{S^0,k}$.
\end{lemma}

\begin{proof}
Since
$$D^j f (s):= i^j \int_0^\infty v^{\alpha+j} e^{ivs}\, dv\ ,$$
it suffices to show that
$$|f(s)|\le C |s|^{-\alpha-1}\ .$$
Assume without loss of generality that $s>0$.
By a change of variables, we need to show
$$\int_0^\infty g\left(\frac {v}{s}\right) v^\alpha e^{iv}\, dv\le C\ .$$
Let $\chi$ be again a smooth cutoff function which is constant
$1$ on $(-\infty, 1]$ and vanishes on $[2,\infty)$.
It suffices to estimate separately the terms
$$\int_0^\infty \chi(v) g\left(\frac {v}{s}\right) v^\alpha e^{iv}\, dv \ ,$$
$$\int_0^\infty (1-\chi(v)) g\left(\frac {v}{s}\right) v^\alpha e^{iv}\,
dv\ .$$
The first term is clearly bounded.
The second term, after $N$ integrations by part, can be estimated by
$$\int_0^{2s}
\left| D^N\left[(1-\chi(v)) g\left(\frac {v}{s}\right) v^\alpha
\right]\right|\, dv\ .$$
By Leibniz'  rule, this is dominated by a constant times
$$\int^{2s}_1 v^{\alpha-N}dv\, ,$$
which is finite if $N$ is choosen sufficiently large.
\end{proof}
This completes the proof of Proposition \ref{r1}

\begin{proposition}\label{o<Rle1}
Assume that $0<R\le 1.$
\be
\item[(a)] If $n=1$, then
$$F_R(s)=e^{iRs}R^{-\frac{1}{2}}b_{-\frac 12}(s)\, ,$$
where $b_{-\frac 12}\in S^{-\frac{1}{2}}$ uniformly in $R\in (0, 1]$.
\item[(b)] If $n\ge 2$, then
$$F_R(s)=e^{iRs}[R^{-\frac{n}{2}}b_{\frac n2 -1}(s)
+R^{1-n}b_0(s)
]\, ,$$
where
$b_{\frac n2 -1}\in S^{\frac{n}{2}-1}$
and
$b_0\in S^0$ uniformly in
$R\in(0,1]$.
\ee
\end{proposition}

For the proof, we  decompose $F_R=F^1_R+F^2_R$, with
\begin{eqnarray*}
F^1_R (s) & := & \int^{\infty}_R \chi (v) D^l_{\sinh
,v}[e^{isv}](\cosh v-\cosh R)^{-\frac{n}{2}+l}dv\, ,\\ F^2_R (s) & := &
\int^{\infty}_2(1-\chi (v))D^l_{\sinh ,v}[e^{isv}](\cosh v-\cosh
R)^{-\frac{n}{2}+l}dv\,
.
\end{eqnarray*}
Here, $\chi$ is a smooth cut-off function such that
$$\chi (v) =1,\mbox{if  }\, \mid v\mid\le 2,\ \text{ and}\ \chi (v) =
0,\,\mbox{if }\mid v\mid\ge 4\, .$$

\medskip

The function $F^2_R$ can be again estimated by
means of Lemma \ref{dshiv} as in Lemma \ref{globg}, which shows that
$F_R^2$ is
in $\S(\R),$ uniformly in $R.$
Thus it remains to estimate $F^1_R$.

\begin{lemma}\label{Dlshiv}
For $0<v\le 4$, we can write
$$D^l_{\sinh ,v}[e^{isv}]=\sum^l_{k=0} s^k q_k (v)v^{k-2l} e^{isv}\ ,$$
where $q_k\in S^0$.
\end{lemma}

\begin{proof} This follows by induction, the case $l=0$ being trivial.
Assume the statement is true for some $l\in\mathbb N$. Then
$$ D^{l+1}_{\sinh,v} [e^{isv}]= \sum^l_{k=0}\, s^k\, D_v
( (\sinh v)^{-1} \, {q_k(v)}\, v^{k-2l} \,e^{isv}). $$
However, on the interval $[0,4]$,
$$\sinh(v)^{-1}=g(v) v^{-1}$$
for some $g\in S^0$. This easily implies the statement of the lemma for $l+1$.
\end{proof}

\begin{lemma}\label{cvcr}
For $0\le v\le 4$, we have
\begin{equation}\label{chv-chR}
 \cosh v-\cosh R  =   \gam(v)(v+R)(v-R)
\end{equation}
for some $\gamma\in S^0$ with $\gamma(v)>\ve>0$
for $0\le v\le 4$, all uniformly in $R\in (0,1]$.
\end{lemma}

\begin{proof}
This follows immediately from
$$\cosh v-\cosh R  = \sum_{n=1}^\infty \frac{(v^2)^n-(R^2)^n}{(2n)!}.$$
\end{proof}

We can therefore decompose
$$F^1_R (s)={\dst\sum^l_{k=0}}
s^k{\dst\int^{\infty}_R}\tilde{\gam}_{k,R}(v)
v^{k-2l}(v+R)^{-\frac{n}{2}+l}(v-R)^{-\frac{n}{2}+l}e^{isv}dv$$
\begin{equation}\label{gkr}
={\dst\sum^l_{k=0}}
s^k \, e^{iRs} {\dst\int^{\infty}_0}\gam_{k,R}(v)
(R+v)^{k-2l}(2R+v)^{-\frac{n}{2}+l}v^{-\frac{n}{2}+l}e^{isv}dv\, ,
\end{equation}
with $\tilde{\gam}_{k,R}, \gam_{k,R}\in S^0$ uniformly in $R\in (0,1]$
and $\gam_{k,R}(v)=0$ for $v>4$.

\begin{lemma}\label{deltalemma}
Let $\al_1,\al_2,\beta\in\R$ such that $\alpha=\al_1+\al_2 >0$ and $-1<\beta $.
Let $\gam\in S^0$ such that $\gam(v)=0$ for $v>4$.
Consider the function
$$f(s)= \int_0^\infty \gam(v) (R+v)^{-\al_1}(2R+v)^{-\al_2} v^\beta
e^{isv}\, dv\ .$$
Then, for every $\delta\ge 0$ with
\begin{equation}\label{delta}
\delta <\alpha\;\mbox{    and    }\;\delta\le\beta +1\, ,
\end{equation}
we have that $R^{\alpha-\delta}f$ is in $S^{\delta-\beta-1}$ uniformly in
$R\in (0,1]$.

\end{lemma}

\begin{proof}
Similarly as in the proof of Lemma \ref{gam}, it suffices to prove that
\begin{equation}\label{fRs}
\mid f(s)\mid\le C
R^{-\alpha+\delta}(1+\mid s\mid)^{\delta-\beta-1},\quad
s\in{\mathbb R}\, ,
\end{equation}
with $C$ independent of $R$.

If $\delta=\beta+1$, then $-\alpha+\beta<-1$ and we have
$$|f(s)|\lesssim
R^{-\alpha+\beta+1} \int_0^\infty (1+v)^{-\alpha} v^\beta \, dv\ ,$$
which proves the desired estimate.

Now assume $\delta<\beta+1$.
Consider first the case $|s|\le 1$.
We estimate
\begin{equation}\label{alphadelta}
(R+v)^{-\al_1}(2R+v)^{-\al_2}\lesssim(R+v)^{-\al} =
(R+v)^{\delta-\alpha}(R+v)^{-\delta}\le R^{\delta-\alpha}v^{-\delta}.
\end{equation}
 This gives
$$|f(s)|\lesssim \int_0^4 R^{\delta-\alpha} v^{\beta-\delta} \, dv\ ,$$
which proves the desired estimate in view of $\delta<\beta+1$.

It remains to consider $|s|>1$.
We write $f=f_1+f_2$ with
$$f_1(s)= s^{-\beta-1}\int_0^\infty \chi(v) \gam(v/s)
(R+v/s)^{-\al_1}(2R+v/s)^{-\al_2}
v^\beta e^{iv}\, dv\ \ ,$$
$$f_2(s)= s^{-\beta-1}\int_1^\infty (1-\chi(v)) \gam(v/s)
(R+v/s)^{-\al_1}(2R+v/s)^{-\al_2}
v^\beta e^{iv}\, dv\ \ ,$$
where $\chi$ is a smooth cutoff function which is constant $1$ on
$(-\infty, 1]$ and vanishes on $[2,\infty)$.
To estimate $f_1$, we split $(R+v/s)^{-\alpha}$ analogously
to (\ref{alphadelta}) and obtain
$$|f_1(s)|\lesssim
s^{-\beta-1}\int_0^\infty |\chi(v)| s^{\delta}R^{\delta-\alpha}
v^{\beta-\delta} \, dv\ \ .$$
This proves the desired estimate for $f_1$ in view of $\delta< \beta+1$.

To estimate $f_2$, we do $N$ times partial integration.
The functions $(1-\chi(v))\gam(v/s)$ is in $S^0$
uniformly in $|s|>1,$ and $v^\beta$ is in $S^\beta$.
With
$$|D_v^k[(R+v/s)^{-\al_1}(2R+v/s)^{-\al_2}]|
\lesssim  (v/s)^k (R+v/s)^{-\alpha-k}v^{-k}$$
$$\lesssim
(R+v/s)^{-\alpha}v^{-k}\lesssim s^{\delta} R^{\delta-\alpha}
v^{-\delta-k}$$
we therefore obtain
$$|f_2(s)|\lesssim s^{\delta-\beta-1} R^{\delta-\alpha} \int_1^\infty
v^{\beta-N}\, dv\ ,$$
which proves the desired estimate for $f_2$.

\end{proof}

We are now in a position to estimate $F_R^1$.
Assume first that $n$ is even and choose $l=n/2$.

Applying Lemma \ref{deltalemma} to (\ref{gkr}) with $\delta=0$ in case $k=l$
and $\delta=1$ in case $k<l$ gives
$$F_R^1=\left[\sum_{k=1}^{l-1} s^k e^{iRs} R^{k-l-\frac n2+1}b_{\frac
n2-l}\right]
+ s^l e^{iRs} R^{-\frac n2} b_{\frac n2 - l -1}$$
$$=\left[\sum_{k=1}^{l-1} e^{iRs} R^{k-l-\frac n2+1}b_{k+\frac n2-l}\right]
+ e^{iRs} R^{-\frac n2} b_{\frac n2  -1}\ ,$$
where $b_\alpha$ generally denotes a function in $S^\alpha$
(uniformly in $R\in (0,1]$) which
may be different at different places in the argument.
Applying Lemma \ref{movers} gives
$$F_R^1=e^{iRs}[ R^{-l-\frac n2+1}b_{\frac n2-l}
+ R^{-\frac n2} b_{\frac n2 -1}]\ .$$
As $l=n/2$, we obtain the desired estimate.

Assume next that $n$ is odd and $n\ge 3$. We choose $l=(n-1)/2$.
Then we apply Lemma \ref{deltalemma} with $\delta=0$ for $k=l$ and
with $\delta=1/2$ for $k<l$ and obtain
$$F_R^1=\left[\sum_{k=1}^{l-1} s^k e^{iRs} R^{k-l-\frac {n-1}2}
b_{\frac {n-1}2-l}\right]
+ s^l e^{iRs} R^{-\frac {n}2} b_{\frac n2 - l -1}$$
$$=\left[\sum_{k=1}^{l-1} e^{iRs} R^{k-l-\frac {n-1}2}b_{k+\frac
{n-1}2-l}\right]
+ e^{iRs} R^{-\frac {n}2} b_{\frac n2 -1}\ .$$
Applying Lemma \ref{movers} again and using $l=(n-1)/2$ gives the
desired estimate.

Finally, assume $n=1$. We choose $l=0$. Applying Lemma
\ref{deltalemma} with $\delta=0$ gives
$$F_R^1=
 e^{iRs} R^{-\frac {n}2} b_{\frac n2 -1}\ ,$$
which proves the desired estimate.

\section{Spectrally localized estimates for the wave propagator}\label{specloc}

The following theorem states pointwise estimates for the convolution
kernel of spectrally localized wave propagators on the $ax+b$ group.

\begin{theorem}\label{ktsatz}
Let $t\in \mathbb R, \lambda > 0$, and let
$\psi$ be an even  bump function $\psi \in C_0^{\infty}(\mathbb
R)$ supported in $[-2,2]$. If $\lambda\ge 1$, we shall in addition
assume that $\psi$ vanishes on $[-1,1]$.
Then the convolution kernel of
$$m_\lambda^t(L):=\psi(\frac{\sqrt{L}}\lambda) \cos(t\sqrt{L})$$
is of the form
$$k_\lambda^t(x,y)=
 e^{-\frac{nx}{2}} e^{-\frac{nR}2}
\left[G_\lambda(R,R-t)+G_\lambda(R,R+t)\right]\ ,$$
where the function $G_\la$ satisfies for every $N\in \N$ the following
estimates:
\begin{enumerate}
\item[(a)]
 If $R\geq 1$, then
$$|G_{\lambda}(R,\rho)|
\lesssim
\begin{cases}
\lambda^{\frac{n}{2}+1} (1+|\lambda \rho|)^{-N},& \text{if
}\lambda \geq 1,\\
\lambda^{2} (1+|\lambda \rho|)^{-N},& \text{if }\lambda < 1.\end{cases}$$
\item[(b)] If $0\leq R \leq 1$, then, for $n=1$,
$$|G_{\lambda}(R,\rho)| \lesssim \begin{cases}R^{-\frac{1}{2}}
\lambda^{\frac{3}{2}} (1+|\lambda \rho|)^{-N},& \text{if }\lambda \geq 1,\\
R^{-\frac{1}{2}} \lambda^{2} (1+|\lambda \rho|)^{-N},& \text{if }
\lambda < 1,\end{cases}$$ and for $n\geq 2$,
$$|G_{\lambda}(R,\rho)| \lesssim \begin{cases}(R^{1-n} \lambda^2 +
 R^{-\frac{n}{2}} \lambda^{\frac{n}{2}+1}) (1+|\lambda \rho|)^{-N},& \text{if
}\lambda \geq 1,\\ R^{1-n} \lambda^{2} (1+|\lambda \rho|)^{-N},& \text{if
}\lambda < 1,\end{cases}$$
\end{enumerate}
where the constants in these estimates depend only on the $C^N$- norms of
$\psi.$
\end{theorem}

\begin{proof}
We consider first the case $R\ge 1$.
By Proposition \ref{nullunend} and
Proposition \ref{r1}, the kernel $k_\lambda^t$ can be written as
$$\frac{c_l}2\, e^{-\frac{nx}{2}}\, e^{-\frac {nR}{2}}\,
\int_\R\, \psi(\frac s\lambda)\, b_{\frac n2 -1}(s)\, s\,
[e^{i(R-t)s}+e^{i(R+t)s}]\, ds\ ,
$$
where $b_{\frac n2 -1}$ is in $S^{\frac n2 -1}$ uniformly in $R\ge 1$.
Then the desired estimate follows by an application of Lemma
\ref{betainS} below with $j=2$.

The case $0<R\le 1$ is done similarly using Proposition \ref{o<Rle1}
instead of Proposition \ref{r1}.
\end{proof}

\begin{lemma}\label{betainS} Let $b\in S^{\beta}$, let
$j\ge 1$ be an integer, and let $\lambda>0$.
Consider
$$ M_{\lambda}(\rho) := \int \psi(\frac{s}{\lambda}) b(s)
s^{j-1} e^{i\rho s} ds, \quad \rho \in \mathbb R.$$ Then, for every $N\in
\mathbb
N$,
$$ |M_{\lambda}(\rho)| \leq C_N \begin{cases} \lambda^{\beta+j}(1+|\lambda
\rho|)^{-N}, & \text{if } \lambda \geq 1 \ \text{and}\
\supp\psi\subset[-2,-1]\cup[1,2],\\
\lambda^{j} (1+|\lambda
\rho|)^{-N}, & \text{if } \lambda<1 \ \text{and}\
\supp\psi\subset[-2,2].\end{cases} $$ Here, the constants $C_N$
depend only on
$N$ and the semi-norms $||b||_{S^{\beta},k}$ of $b$, and on the $C^N$-norm of
$\psi$.
\end{lemma}

\begin{proof} In order to defray the notation, we shall write $A \lesssim
B$, if
$A\leq C_N B$, where $C_N$ is an ''admissible'' constant in the
sense  described in the lemma. We may and
shall assume that $\rho \geq 0$. We write
$$M_{\lambda}(\rho) = \lambda^{j} \int \psi(s) \,
b(\lambda s) \,s^{j-1} \,e^{i\la \rho s}\, ds\ .$$

Assume $\lambda \ge 1$ and $\psi$ is supported
in $[-2,-1]\cup [1,2]$.
By the symbol estimates for $b$ we have for $s$ in the support of $\psi$:
$$|D^k_s b(\lambda s)|
\lesssim \lambda^{k} (1+|\lambda s|)^{\beta-k}
\lesssim \lambda^{\beta} (\frac 1\lambda +|s|)^{\beta-k}\lesssim
\lambda^{\beta}\ .$$
Integrating by parts $N$ times, we thus find that
$|M_{\lambda}(\rho)| \lesssim \lambda^{\beta+j}(1+ |\lambda \rho|)^{-N}$.

Next, assume $\lambda \le 1$ and $\psi$ is supported
in $[-2,2]$.
We then have for $s$ in the support of $\psi$
$$|D^k_s b(\lambda s)|
\lesssim \lambda^{k} (1+|\lambda s|)^{\beta-k}\lesssim 1\ \ \ .$$
Integrating again by parts $N$ times, we find that
$|M_{\lambda}(\rho)| \lesssim \lambda^{j}(1+ |\lambda \rho|)^{-N}$.
\end{proof}

 As a consequence of Theorem \ref{ktsatz}, we
obtain estimates of the $L^1$-norms of the convolution  kernels of
$\psi(\frac{\sqrt{L}}{\lambda}) \cos(t\frac{\sqrt{L}}{\lambda}).$ Notice that
the corresponding multipliers  result from a re-scaling  of the multiplier
for the case
$\lambda=1,$ but the kernels cannot just be obtained by some scaling argument
from the case $\lambda=1,$  since
the operator $L$ is not homogeneous.

\begin{proposition}\label{wlambdat} Let $W^t_{\lambda} :=
k^{t/\lambda}_{\lambda}$ denote the convolution kernel of
 $\psi(\frac{\sqrt{L}}{\lambda})
\cos(t \frac{\sqrt{L}}{\lambda})$, and let $\ve\geq 0$.
\begin{itemize}
\item [(a)] If $\lambda \geq 1$ and $\supp \psi\subset [-2,-1]\cup [1,2],$
then

\begin{equation} \label{wlamb1}
\int \limits_G |W_{\lambda}^t (x,y)| R(x,y)^{\ve} d(x,y) \lesssim
\begin{cases}
\lambda^{-\ve} (1+t)^{\frac{n}{2} + \ve}, & \text{if } t\leq \lambda,\\
\lambda^{\frac{n}{2}-1-\ve} ~t^{1+\ve}, & \text{if } t\geq \lambda.
\end{cases}
\end{equation}

\item [(b)] If $0<\lambda\leq 1$ and $\supp \psi\subset [-2,2],$
then
\begin{equation}\label{wlamb2}
\int\limits_G |W^t_{\lambda}(x,y)| R(x,y)^{\ve} d(x,y) \lesssim \lambda^{-\ve}
(1+t)^{1+\ve}.
\end{equation}
\end{itemize}
In particular, if $\lambda \geq 1$, then
\begin{equation}\label{wlamb3}
\int\limits_G |W^t_{\lambda}(x,y)|~ (1+ \lambda R(x,y))^{\ve} ~d(x,y)
\lesssim \begin{cases}(1+t)^{\frac{n}{2}+\ve}, & \text{if } n \geq2,\\
(1+t)^{1+\ve}, & \text{if } n=1,\end{cases}
\end{equation} and,
 if $0<\lambda \leq 1$, then
\begin{equation}\label{wlamb4}
 \int\limits_G
|W^t_{\lambda}(x,y)|~ (1+ \lambda R(x,y))^{\ve} ~d(x,y)
\lesssim (1+t)^{1+\ve},
\end{equation}
 in each instance  uniformly in $\lambda$. The constants in these estimates
depend only on the $C^N$-norms of $\psi.$
\end{proposition}

\begin{proof} Without loss of generality, we shall assume that $t\geq 0$.
Then the
dominant term in Theorem \ref{ktsatz} is the one containing
$G_{\lambda} (R,R-t)$, to which we
shall therefore  restrict ourselves.

We consider first $\lambda\ge 1$.
By Theorem \ref{ktsatz} and Lemma \ref{intlem}, we can estimate the
left-hand-side of (\ref{wlamb1}) by
$$
\lambda^2 \int_0^1
(1+|\lambda R-t|)^{-N} R^{1+\ve}\, dR
+\lambda^{\frac n2 + 1} \int_0^1
(1+|\lambda R-t|)^{-N} R^{\frac n2 + \ve}\, dR$$
\begin{equation}\label{threeterms}
+ \lambda^{\frac n2 + 1} \int_1^\infty
(1+|\lambda R-t|)^{-N} R^{1+\ve}\, dR\ .\end{equation}
If $t\le \frac \lambda 2$, then we can estimate this using Lemma \ref{alphagr0}
below by
$$
\lambda^{-\ve}(1+t)^{1+\ve}+\lambda^{-\ve}(1+t)^{\frac n2 + \ve}
+\lambda^{\frac n2 +1 -N}\ .$$
Similarly, if $\frac \lambda 2 \le t\le 2 \lambda$, we estimate
(\ref{threeterms}) by
$$
\lambda^{-\ve}(1+t)^{1+\ve}+\lambda^{-\ve}(1+t)^{\frac n2 + \ve}
+\lambda^{\frac n2 -1 -\ve}(1+t)^{1+\ve}\ ,$$
and if $2\lambda \le t$ we estimate (\ref{threeterms}) by
$$
(1+t)^{-N}
+\lambda^{\frac n2 -1 -\ve}(1+t)^{1+\ve}\ .$$
In each case we easily verify (\ref{wlamb1}), taking into account that
for $n=1$ the first of the three summands of (\ref{threeterms}) is
not present.

If $0<\lambda\le 1$, we estimate the
left-hand-side of (\ref{wlamb2}) by
$$
\lambda^2 \int_0^\infty
(1+|\lambda R-t|)^{-N} R^{1+\ve}\, dR$$
if $n\ge 2$ and by
$$
\lambda^2 \int_0^1
(1+|\lambda R-t|)^{-N} R^{\frac 12+\ve}\, dR
+ \lambda^2 \int_1^\infty
(1+|\lambda R-t|)^{-N} R^{1+\ve}\, dR$$
if $n=1$. In either case it is easy to verify
(\ref{wlamb2}) using Lemma \ref{alphagr0}.
Estimates (\ref{wlamb3}) and (\ref{wlamb4}) follow
immediately from estimates (\ref{wlamb1}) and (\ref{wlamb2}).
\end{proof}

\begin{lemma}\label{alphagr0}
For $\alpha \geq 0$, $t \geq 0$ and $N>\alpha + 1$,
let
\begin{eqnarray*}I_0 &:=& \int \limits^1_0 (1+ |\lambda R-t|)^{-N}
R^{\alpha}\, dR,\\ I_{\infty} &:=& \int \limits_1^{\infty}(1+ |\lambda
R-t|)^{-N} R^{\alpha}\, dR.
\end{eqnarray*}
 Then
\begin{equation}\label{inull}
I_0 \leq C
\begin{cases}
 (1+\lambda)^{-\alpha-1} (1+ t)^{\alpha},
&\text{if } t \leq 2\lambda,
\\(1+ t)^{-N}, & \text {if } t \geq 2\lambda,
\end{cases}
\end{equation} and
\begin{equation}\label{iunend} I_{\infty} \leq
C\begin{cases}\lambda^{-\alpha-1}(1+\lambda)^{-N+\alpha+1}, & \text{if } t\leq
\frac{\lambda}{2},\\
 \lambda^{-\alpha-1}(1+ t)^{\alpha}, & \text{if }
 t\geq \frac{\lambda}{2}.
\end{cases}
\end{equation}
\end{lemma}

\begin{proof}
 We begin with $I_0$.
 If $t\geq 2\lambda$, then clearly $I_0 \lesssim (1+
t)^{-N}$.
 If $t \leq 2\lambda$, we write
$$I_0 = \lambda^{-\alpha-1} \int \limits ^{\lambda -  t}_{- t}
(1+ |v|)^{-N} (v+  t)^{\alpha} dv.$$  If $\lambda \geq 1$,
 one easily deduces from this representation that $I_0 \lesssim
\lambda^{-\alpha-1} (1+  t)^{\alpha}$. And, if $\lambda\leq 1$, the
original formula for $I_0$
 immediately implies $I_0 \lesssim 1$, so that we obtain (\ref{inull}).

 As for $I_{\infty}$, if $t\leq \frac{\lambda }{2}$, then clearly
$$I_{\infty} \lesssim \int \limits_{1}^{\infty} (1+\lambda R)^{-N} R^{\alpha}\,
dR = \lambda^{-\alpha-1} \int \limits_{\lambda}^{\infty} (1+R)^{-N}
R^{-\alpha}\, dR,$$ hence $I_{\infty} \lesssim \lambda^{-\alpha-1}
(1+\lambda)^{-N+\alpha+1}$.\\ If $t\geq \frac{\lambda }{2}$, we write
$$I_{\infty} = \lambda^{-\alpha-1} \int \limits^{\infty}_{\lambda- t}
(1+|v|)^{-N} (v+ t)^{\alpha} dv.$$
If $ t \leq 1$, this implies $I_{\infty} \lesssim
\lambda^{-\alpha-1}$, and if $ t \geq 1$, one finds that
$I_{\infty} \lesssim \lambda^{-\alpha-1} t^{\alpha}$, so
that also (\ref{iunend}) is verified.\par
\end{proof}

By means of the subordination principle described e.g. in \cite{muller}, we
immediately obtain:

\begin{corollary}\label{hs61}(cf. \cite{hebisch-steger},Theorem. 6.1)\\
If $\ve>0$, $s_0, s_1>\frac{3}{2}+\ve$ and $s_1 > \frac{n+1}{2} + \ve$,
then there exists a constant $C$ such that, for every continuous function
$F$ supported in $[1,2]$ and $0<\lambda \leq 1,$ $$ \int\limits_G
|F(\frac{L}{\lambda^2}) \delta_0(x,y) |
 ~(1+\lambda R(x,y))^{\ve} ~d(x,y) \leq C ||F||_{H(s_0)},$$
while for $\lambda \geq 1$
$$ \int\limits_G |F(\frac{L}{\lambda^2}) \delta_0(x,y) | ~(1+\lambda
R(x,y))^{\ve}~ d(x,y) \leq C ||F||_{H(s_1)}.$$
\end{corollary}

\begin{proof} Choose an even function $\psi \in C_0^{\infty}(\R)$ such that
$\psi=1$ on
$[1,2]$ and
$\supp\psi\subset [-4,-\frac12]\cup[\frac{1}{2},4]$. Proposition \ref{wlambdat}
holds for such
$\psi$ as well. Put $f(v):= F(v^2)$. Then $||f||_{H(s)} \sim ||F||_{H(s)}$, for
any
$s\geq 0$, and $F(\frac{L}{\lambda^2}) = f(\frac{\sqrt{L}}{\lambda})=
\psi(\frac{\sqrt{L}}{\lambda}) f(\frac{\sqrt{L}}{\lambda})$. Moreover,
by the Fourier inversion formula and Fubini's theorem, one easily obtains
$$f(\frac{\sqrt{L}}{\lambda}) = \frac{1}{\pi} \int \limits_0^{\infty}
\hat{f}(t)
\cos(t \frac{\sqrt{L}}{\lambda}) \,dt,$$
since $f$ is an even function. Thus $$F(\frac{L}{\lambda^2}) =
\frac{1}{\pi} \int\limits_0^{\infty}  \hat{f}(t)
\psi(\frac{\sqrt{L}}{\lambda}) \cos(t \frac{\sqrt{L}}{\lambda})\, dt,$$
which implies \begin{eqnarray*}I_{\lambda}& := & \int\limits_G
|F(\frac{L}{\lambda^2})
\delta_0 | (1+ \lambda R)^{\ve} d(x,y)\\
&\lesssim& \int\limits_0^{\infty} |\hat{f}(t)| \left[\int |W_{\lambda}^t|
(1+\lambda R)^{\ve} d(x,y)\right] dt.\end{eqnarray*} Thus, if $0<\lambda \leq
1$, then, by (\ref{wlamb4}),

\begin{eqnarray*} I_{\lambda} & \lesssim & \int \limits_{\mathbb R}
|\hat{f}(t)| (1+|t|)^{1+\ve} dt\\
&\lesssim& \left( \int\limits_{\mathbb R} \left|\hat{f}(t)
(1+|t|)^{s_0}\right|^2 dt\right)^{\frac{1}{2}}\\ &=& ||f||_{H(s_0)} \sim
||F||_{H(s_0)}.
\end{eqnarray*}The case $\lambda \geq 1$ can be treated in the same way.
\end{proof}

For the class of groups G considered here, we have thus established a
completely
 different approach to the basic Theorem 6.1 in \cite{hebisch-steger}, entirely
based on the wave equation.\\

\section{ Improvements on the estimates in Theorem \ref{ktsatz} for small $R$}

The estimates in Theorem \ref{ktsatz} are already good enough for
$L^1$-estimates,
but not yet for $L^\infty$-estimates, since they exhibit singularities at
$R=0.$
One knows that the singular support of the wave propagator for time $t$  is
the
sphere $R=|t|,$ so that these singularities are in fact  not present. We
shall show in
this section how to improve on our estimates
when $R\le 1,$ which we shall assume throughout this section.

To this end we observe that by formula \eqref{sub}, we may  replace
$F_R(s)$ in the previous
discussions by $F_R(s)-F_R(-s)=2i\tilde F_R(s),$ where
$$\tilde F_R(s):=\int\limits^{\infty}_{R}
D^l_{\sinh, v}[\sin(sv)](\cosh v-\cosh R)^{-\frac{n}{2}+l}dv,\quad (s\ge 0,
\quad l>\frac{n}{2}-1). $$

Working with $\tilde F_R(s)$ in place of $ F_R(s),$ we can prove the following
theorem.

\begin{theorem}\label{2ktsatz}
If $R\le 1,$ then the estimates in Theorem \ref{ktsatz} can be improved by
the following
 additional estimates, valid for any $N\in\N :$
\begin{equation}\label{newest}
|G_{\lambda}(R,\rho)| \lesssim \begin{cases}
\lambda^{n+1} (1+|\lambda \rho|)^{-N},& \text{if }\lambda \geq 1,\\
 \lambda^{2} (1+|\lambda \rho|)^{-N},& \text{if }
\lambda < 1,\end{cases}
\end{equation}
\end{theorem}

\medskip
In order to prove this result, as in Section  \ref{refine} we
 split $\tilde F_R(s)=\tilde F^1_R(s)+
\tilde F^2_R(s),$ where
$$\tilde F^1_R(s):=\int\limits^{\infty}_{R}\chi(v)
D^l_{\sinh, v}[\sin(sv)](\cosh v-\cosh R)^{-\frac{n}{2}+l}dv. $$

$\tilde F^2_R(s)$ is again a Schwartz function, uniformly for $0\le R\le 1,$
and its contribution to $k_\la^t$ can easily be seen to satisfy the estimates
in Theorem \ref{2ktsatz}.

In order to deal with $\tilde F^1_R(s),$ we need the following substitute for
Lemma \ref{Dlshiv}.

\begin{lemma}\label{2Dlshiv}
For $0<sv<\pi/2$ and $0<v<4$ we can write
\begin{equation}\label{sh2}
D^l_{\sinh, v}[\sin(sv)]=\sum_{k=0}^l s^{2k} q_k(sv,v),
\end{equation}
where $q_k(y,v)$ has a power series expansion of the form
\begin{equation}\label{sh3}
\sum_{m,n=0}^\infty a_{mn} y^{2m+1} v^{2n}.
\end{equation}
Moreover,
\begin{equation}\label{sh4}
|D_s^j[q_k(sv,v)]|\le C_{j,k}(1+s^2)^{-j/2}, \quad\mbox{uniformly in }  v.
\end{equation}
\end{lemma}

\begin{proof}We proceed by induction on $l,$ the case $l=0$ being clear.
Assume that $q_k(y,v)$ is given by \eqref{sh3}.  Then
$$\frac {q_k(sv,v)}{\sinh(v)}=s g_k(sv,v),$$
where $g_k(y,v)$ has an expansion of the form
$$g_k(y,v)=\sum_{m,n=0}^\infty b_{mn} y^{2m} v^{2n}.$$
Then
\begin{eqnarray*} D_v[sg_k(sv,v)]=&&\sum_{m\ge 1, n\ge 0} b_{mn}
 2m s^2 (sv)^{2m-1} v^{2n}\\
&&+\sum_{m\ge 0, n\ge 1} s b_{mn} 2n  (sv)^{2m} v^{2n-1}\\
=&& s^2 h_k^1(sv,v) +h^2_k(sv,v),
\end{eqnarray*}
where $h_k^1(y,v)$ and $h_k^2(y,v)$ are of the form \eqref{sh3}. This shows
that
\eqref{sh3} holds also for $l+1$ in place of $l.$

Moreover, \eqref{sh4} is obvious for $|s|\le 1,$ in view of \eqref{sh3}.
And, if
$|s|\ge 1,$ it follows from
$$D_s^j[q_k(sv,v)]=v^j (D_y^j q_k)(sv,v)=s^{-j} (sv)^j(D_y^j q_k)(sv,v).$$

\end{proof}
Now, if  $\la R\ge 1/2$, then $\la\ge 1/2,$ and the estimates \eqref{newest}
follow immediately from Theorem \ref{ktsatz}. Let us therefore assume that
$\la R\le 1/2.$

We write
$$\tilde F^1_R(s)=H^1_R(s)+H^2_R(s),$$
where
$$H^1_R(s):= \int_R^{\frac 1{2\la}}\chi(v) D^l_{\sinh, v}[\sin(sv)]
(\cosh v-\cosh R)^{-\frac{n}{2}+l}dv, $$
$$H^2_R(s):= \int_{\frac 1{2\la}}^\infty \chi(v) D^l_{\sinh, v}[\sin(sv)]
(\cosh v-\cosh R)^{-\frac{n}{2}+l}dv. $$

\noindent We need information on the asymptotics of these functions for
$|s|\le 2\la.$

\medskip
 As for $H^1_R(s),$ notice that in the integral defining $H^1_R(s)$ we have
$|sv|\le 1$ if $|s|\le 2\la.$ Therefore, from Lemma \ref{2Dlshiv} and Lemma
\ref{cvcr}
 we find that,
 for any $l>-n/2+1,$

$$H^1_R(s)=\sum_{k=0}^l s^{2k}\int_R^{\frac 1{2\la}}\gamma_k(sv,v)
(v+R)^{-\frac n2 +l}(v-R)^{-\frac n2 +l}\, dv,$$
where $\gamma_k(sv,v)$ is supported where $0\le v\le \min (2,\frac 1{2\la})$
and satisfies $\gamma_k(-sv,v)=-\gamma_k(sv,v)$ and
$$|D_s^j[\gamma_k(sv,v)]|\le C_{j,k}(1+s^2)^{-j/2}, \quad\mbox{for every  }
j\in\N.
$$

Choose $l$ large enough so that $l-\frac n2\ge 0, $ and let
$$J(\la):=\int_R^{\frac 4{(1+\la)}}(v+R)^{-\frac n2 +l}(v-R)^{-\frac n2
+l}\, dv.
$$
Then clearly
$$J(\la)\le (1+\la)^{-2l+n-1},$$
and we find that, for $|s|\le 2\la,$
$$H^1_R(s)=(1+\la)^{-2l+n-1}\sum_{k=0}^l s^{2k} b_{0,k}(s),$$
where $b_{0,k}$ is an odd function in $S^0,$ uniformly in $R$ and $\la.$
From Lemma \ref{betainS} we therefore
obtain
\begin{eqnarray*}&&\left|  \int_0^\infty \psi(\frac s\la)H^1_r(s)
s\cos (ts)\, ds\right|=\frac 12\left| \int_\R \psi(\frac s\la)H^1_r(s)
s\cos (ts)\, ds\right|\\
&&\qquad \lesssim C_N\sum_{k=0}^l (1+\la)^{-2l+n-1}\la^{2k+2}(1+|\la t|)^{-N},
\end{eqnarray*}
hence
\begin{eqnarray}\label{este1}
&&\left|  \int_0^\infty \psi(\frac s\la)H^1_r(s)
s\cos (ts)\, ds\right|
\lesssim C_N\begin{cases}
\lambda^{n+1} (1+|\lambda t|)^{-N},& \text{if }\lambda \geq 1,\\
 \lambda^{2} (1+|\lambda t|)^{-N},& \text{if }
\lambda < 1,\end{cases}
\end{eqnarray}

\medskip

Next, we consider $H^2_R(s),$  again for $|s|\le 2\la.$

Observe first that $H^2_R\equiv 0,$ unless $\la\ge 1/4.$ In the latter
case, one finds that
$H^2_R(s)$ behaves like $F^1_R(s),$ only with $R$ replaced by $\frac 1\la
\le 4.$
Replacing $R$ by $\frac 1\la$ and $R\pm t$ by $\pm t$  in Theorem
\ref{ktsatz} (b), we therefore find that
$$\left|  \int_0^\infty \psi(\frac s\la)H^1_r(s)
s\cos (ts)\, ds\right|$$
satisfies estimates \eqref{este1} too.

Noticing finally that $1+|\la t|\sim 1+|\la \rho|$ if $\rho=R\pm t,$ since
 $\la R\le 1/2,$ the conclusion of Theorem \ref{2ktsatz} follows.

\begin{cor} Assume that $\la\ge 1$ and $t\ge 0.$ Then
\begin{equation}\label{supest}
||k^t_\la||_\infty\lesssim (1+t^{-\frac n2}) \la^{\frac n2 +1}.
\end{equation}
\end{cor}

\begin{remark} {\rm Notice that, for small times, this estimate agrees with
the one valid for the Laplacian on Euclidean
space $\R^{n+1},$  as is to be expected, since $L$ is elliptic.  However,
for large times, there
appears no
dispersive effect (definitely not for $n=2$, by Hebisch's transfer
principle), so that it seems unlikely that non-trivial Strichartz-type
estimates
will hold for large times.}
\end{remark}

\begin{proof} First we observe that $e^{-\frac{nx}2}e^{-\frac{nR}2}\le 1,$
and equality
holds, if $y=0$ and $x\le 0.$ Therefore,
$$ ||k^t_\la||_\infty\lesssim \sup_{R\ge 0} |G_\la(R,R-t)|.$$
\medskip

If $R\ge 1,$ then, by Theorem \ref{ktsatz},
$$ |G_\la(R,R-t)|\lesssim \la^{\frac n2 +1}.$$

So, assume that $R\le 1.$ Then, by Theorem \ref{2ktsatz},
\begin{equation}\label{gest}
 |G_\la(R,R-t)|\lesssim \la^{n+1}(1+\la|R-t|)^{-N}
\end{equation}
for every $N\in\N.$

If $\la t\le 1,$ this implies
$$ |G_\la(R,R-t)|\lesssim \la^{n+1}\le\la^{\frac n2 +1}t^{-\frac n2}.$$

Assume next that $\la t\ge 1.$

If $R\le t/2,$ then $|R-t|\sim t,$ so that \eqref{gest} implies

$$ |G_\la(R,R-t)|\lesssim \la^{n+1}(\la t)^{-N}$$
for every $N\in\N,$ hence
\begin{equation}\label{rightest}
 |G_\la(R,R-t)|\lesssim \la^{\frac n2+1}t^{-\frac n2}.
\end{equation}

\medskip
If $ R\ge t/2,$ then  for $n\ge 2$
Theorem \ref{ktsatz}
implies \eqref{rightest}, since
$t^{1-n}\la^2\le t^{-\frac n2}\la^{\frac n2+1},$  and \eqref{rightest} is also
valid for $n=1.$

\end{proof}

\begin{remark}{\rm The group $G$ can be considered as an Iwasawa
$AN$-subgroup of the
Lorentz group $S=SO(1,n+1),$ and hence may be identified as a manifold with
the symmetric space
$K\backslash S,$ where $K$ is a maximal compact subgroup of  $S.$ The
spherical
function $\varphi_0$ of order zero on $K\backslash S$ is comparable to
$(\frac R{\sinh(R)})^{n/2}$ in these coordinates, as one finds from
Harish-Chandra's spherical
 function expansion  (see e.g. \cite{helgason}). In view of the well-known
estimates for the wave propagators in Euclidean space, a  naive extrapolation
of Hebisch's transfer principle to this situation (where $S$ is not a
complex semisimple Lie group,
unless $n=2$) would lead to the following``conjecture'':
$$k_\lambda^t(x,y)=
 e^{-\frac{nx}{2}} e^{-\frac{nR}2} P_\la^t(R),$$
where
\begin{enumerate}
\item[(a)]
 If $R\geq 1$, then
$$|P_\la^t(R)|
\lesssim
\begin{cases}
t^{-\frac n2}\lambda^{\frac{n}{2}+1}R^{\frac n2} (1+\lambda |R-t|)^{-N},&
\text{if
}\lambda t \geq 1,\\
\lambda^{n+1} R^{\frac n2}(1+\lambda R)^{-N},& \text{if }\lambda t<
1.\end{cases}$$
\item[(b)] If $0\leq R \leq 1$, then
$$|P_\la^t(R)| \lesssim \begin{cases}t^{-\frac n2}\lambda^{\frac{n}{2}+1}
(1+\lambda |R-t|)^{-N},& \text{if }\lambda t \geq 1,\\
\lambda^{n+1} (1+\lambda R)^{-N},& \text{if }\lambda t< 1,\end{cases}$$
\end{enumerate}

From Theorems \ref{ktsatz} and \ref{2ktsatz}, one can indeed easily verify
these estimates, if
$\la\ge2,$ say.

However, if $\la\le 1,$ and if we choose $R=t\ge1$ and $\la t\ge 1,$ the
``conjecture'' would predict
a size of order $\la^{\frac n2+1}$  for $|P_\la^t(R)|,$ whereas we find the
order $\la^2.$

}
\end{remark}

\section{Growth estimates for
solutions to the wave equation
in terms of spectral Sobolev norms}\label{sobolev}

\begin{theorem}\label{symbol}
Given a symbol $m \in S^{-\alpha}$, we define operators $T_{1}^t := m(\sqrt{L})
\cos(t\sqrt{L})$ and $T_{2}^t := m(\sqrt{L})
\frac{\sin(t\sqrt{L})}{\sqrt{L}}$, a
priori on $L^2(G)$, for $t\in \mathbb R$. Let $1\leq p\leq
\infty$.
\begin{itemize}
\item[(a)] If $\alpha>n |\frac{1}{p}-\frac{1}{2}|$, then $T_1^t$ extends from
$L^p \cap L^2(G)$ to a bounded operator on $L^p(G)$, and
$$ ||T^t_1||_{L^p\rightarrow L^p} ~\leq~
C_p~(1+|t|)^{2|\frac{1}{p}-\frac{1}{2}|}.$$
\item[(b)] If $\alpha>n|\frac{1}{p}-\frac{1}{2}|-1$, then $T_2^t$ extends from
$L^p\cap L^2(G)$ to a bounded operator on $L^p(G)$, and $$||T^t_2||_{L^p
\rightarrow L^p} ~\leq~ C_p~(1+|t|).$$
\end{itemize}
\end{theorem}

Note that the extension is unique, if $1\leq p<\infty$.\\

\proof (a) Let $\chi \in C_0^{\infty} (\mathbb R)$ be an even function such
that
$\chi(s) =1$ if $|s| \leq \frac{1}{2},$ and $\chi(s)=0$, if $|s|\geq 1$. Put
$\psi_0(s) := \chi(\frac{s}{2})$, and $\psi_j(s) :=\chi(2^{-j-1}s) -
\chi(2^{-j}s)
= \psi(2^{-j}s),$ $ j=1,\ldots\infty,$ where $\psi(s) := \chi(\frac{s}{2}) -
\chi(s)$
 is supported in $\{s:~\frac{1}{2} \leq |s| \leq 2\}$. Then $\psi_0$ is
supported in $[-2,2]$, $\psi_j$ in $\{s:~ 2^{j-1} \leq |s| \leq 2^{j+1} \}$
for $j\geq 1$, and
\begin{equation}\label{sumpsi} \sum\limits_{\j=0}^{\infty} \psi_j(s) =1,
\quad s\in \mathbb R. \end{equation} We shall restrict ourselves to the
case $1\leq p <2$,
since the case $p=2$ is trivial and the case $p>2$ follows from the case $p<2$
 by duality. Using \eqref{sumpsi}, we decompose the symbol $m$ as
$$m(s) = \sum_{j=0}^{\infty} m_j(2^{-j}s),$$
 where $m_0 = m\chi$ and $m_j(s) :=
(m\psi_j)(2^js) = m(2^js) \psi,  \text{ if } j\geq 1.$ Notice that
\begin{equation}\label{mjott} ||~m_j~||_{C^N} ~\leq ~  C~ 2^{-\alpha
j},\end{equation}where the constant $C$ depends on the semi-norms
$||m||_{S^{-\alpha},k}$ only.\\

Then, for every $f \in L^2(G)$, \begin{equation}\label{T1t} T_1^t f ~=~
\sum\limits^{\infty}_{j=0} T_j f\quad \text{in   }L^2(G),\end{equation}
where $T_j := m_j(\frac{\sqrt{L}}{2^j}) \cos ((2^jt) \frac{\sqrt{L}}{2^j}).$\\
Estimating the operator norms $||T_j||_{L^1 \rightarrow L^1}$ of $T_j$ on
$L^1(G)$ by means of Proposition \ref{wlambdat} and (\ref{mjott}), and
interpolating these estimates with the trivial $L^2$-estimate $||T_j||_{L^2
\rightarrow L^2}
\lesssim 2^{-\alpha j}$, we obtain the following inequalities (we assume
w.l.o.g.
$t\geq 0$):
\begin{equation}\label{tjott}||~T_j~||_{L^p \rightarrow L^p} \lesssim~
\begin{cases}
2^{-\alpha j} ~(1+2^jt)^{n|\frac{1}{p}-\frac{1}{2}|}, & \text{if  }t\leq 1,\\
2^{-\alpha j}~(2^{\frac{n}{2}j}t) ^{2|\frac{1}{p}-\frac{1}{2}|}, & \text{if
}t\geq 1. \end{cases}\end{equation} The estimate in (a) follows immediately
from
\eqref{tjott} by summation over all $j\geq 0$.\\

As for (b), observe first that if we replace $m^t_{\lambda}(s) =
\psi(\frac{\sqrt{s}}{\lambda}) \cos(t\sqrt{s})$ in Section \ref{specloc} by
$\tilde{m}^t_{\lambda}(s) = \psi(\frac{\sqrt{s}}{\lambda})
\frac{\sin(t\sqrt{s})}{\sqrt{s}},$ then the factor $s(e^{i(R-t)s} +
e^{i(R+t)s})$ in the corresponding kernel $k^t_{\lambda}(R)$ has to be replaced
by $i(e^{i(R-t)s} - e^{i(R+t)s})$. By Lemma \ref{betainS}, with $j=1$, the
estimates for the function $\tilde{k}^t_{\lambda}$ associated to
$\tilde{m}^t_j$ are therefore the same as for $k^t_{\lambda}$, except for
an additional factor $\lambda^{-1}$. Moreover,
$$\sup\limits_{s} \left|m_j(\frac{s}{2^j})\frac{\sin(ts)}{s}\right|~ \lesssim
~\begin{cases}2^{-\alpha j}~ 2^{-j},& \text{if  } j\geq 1,\\2^{-\alpha
j}~(1+t),&
\text{if  } j=0.\end{cases}$$

Together, this implies that for $j\geq 1$, the operators $\tilde{T}_j$
arising in the dyadic decomposition of $T_2^t$ satisfy the same estimates
as $T_j$,
except for an additional factor $2^{-j}$. And, for $j=0$,
$$||~\tilde{T}_0~||_{L^2 \rightarrow L^2}~\lesssim ~ (1+t)~,\qquad
||~\tilde{T}_0~||_{L^1 \rightarrow L^1}~\lesssim ~ (1+t), $$
hence $||~\tilde{T}_0~||_{L^p \rightarrow L^p} \lesssim (1+t).$ The
estimates in (b)
thus follow by summing over all $j$.
\par

\hfill{Q.E.D}
\medskip

Let $u=u(t,x)=u_t(x)$ be the solution of the Cauchy problem
\begin{equation}\label{cauchy}
\frac{\partial^2}{\partial t^2} u-(X^2 + \sum\limits_{j=1}^n Y_j^2) u =0,
\quad u_0=f, \quad \frac{\partial}{\partial t} u |_{t=0} =g. \end{equation}

Then, a priori for $f, g \in L^2(G)$,  $u_t$ is given by $u_t= \cos(t\sqrt{L})
f + \frac{\sin(t \sqrt{L})}{\sqrt{L}} g.$  If we define adapted Sobolev norms
$$||~\varphi~||_{L^p_{\alpha}} ~:= ~
||~(1+L)^{\frac{\alpha}{2}}~\varphi~||_{L^p}, \quad \alpha \in \mathbb R,$$ we
therefore immediately obtain from Theorem \ref{symbol} the following

\begin{corollary}\label{alpha0} If $1\leq p< \infty$, then for $\alpha_0 >
n~|\frac{1}{p}-\frac{1}{2}|$ and $\alpha_1>n~|\frac{1}{p}-\frac{1}{2}|-1$
\begin{equation}\label{utlp}
||u_t||_{L^p} ~\leq ~ C_p
(~(1+|t|)^{|\frac{2}{p}~-1|}||f||_{L^p_{\alpha_0}}~+~(1+|t|)~||g||_{L^p_{\alpha_
1}})\ \ .
\end{equation}\end{corollary}

\begin{remark}It is likely that the estimate (\ref{utlp}) even holds for
$\alpha_0=n~|\frac{1}{p}-\frac{1}{2}|$ and
$\alpha_1=n~|\frac{1}{p}-\frac{1}{2}|-1$,
if $1<p<\infty$. This would be the counterpart to corresponding results by
Peral
\cite{peral} and Miyachi \cite{miyachi} in the Euclidean setting
 (see also \cite{sss} for a local variable coefficient version). The sharp
result would require an introduction of a suitable Hardy respectively BMO-space
on
$G$. There is strong evidence that such spaces  exist on $G$, in view of the
ideas in \cite{ionescu} and \cite{hebisch-steger}, but we shall not pursue
these issues here. \end{remark}

\bibliographystyle{plain}
\bibliography{axbJAN04}

\end{document}